\documentclass[12pt,thmsa]{article}
\usepackage{amsfonts}
\usepackage{amssymb}
\newtheorem{teo}{Theorem}[section]

\newtheorem{obs2}[teo]{Remark}

\newtheorem{tea}{Theorem}[subsection]

\newtheorem{no2}[teo]{Note}

\newtheorem{no3}[tea]{Note}

\newcommand{\Frob}{{\rm Frob }}

\newcommand{\mod}{{\rm mod}}

\newcommand{\Q}{\mathbb{Q}}

\newcommand{\PGL}{{\rm PGL}}

\newcommand{\eq}{x^4 + y^4 = q z^p}

\title{Solving diophantine equations  $x^4 + y^4 = q z^p $}
\author{Luis V. Dieulefait  \thanks{supported
by a MECD postdoctoral grant at the Centre de Recerca Matem{\'a}tica
from Ministerio de Educaci{\'o}n y Cultura}\\
Centre de Recerca Matem{\'a}tica \\  Apartat 50, E-08193 Bellaterra, Spain\\
e-mail: LDieulefait@crm.es}

\begin{document}

\maketitle
%\tableofcontents
%\vspace{6mm}
%\newpage
%Author's name and mailing address: see page 1\\
%e-mail: luisd@cerber.mat.ub.es
%\newpage
\begin{abstract}
We give a method to solve generalized Fermat equations of type
$\eq$, for some prime values of $q$ and every prime $p$  bigger
than $13$. We illustrate the method by proving that there are no
solutions for $q= 73, 89$ and $113$.\\
Math. Subject Classification: 11D41,11F11
\end{abstract}
%\newpage

\section{Introduction}
When considering the non-existence of primitive solutions in
integers for the equation:
$$ \eq  \quad \qquad (1.1) $$
for a fixed prime $q$ and $p$ running over all sufficiently large
prime exponents, there are two sets of values of $q$  where the
answer is immediate. On the one hand, for primes $q$ that can be
written as the sum of two biquadrates, like $17 = 1 + 16$, one
already has a trivial solution for the equation above for an
arbitrary value of $p$ just taking $z=1$. On the other hand, a
necessary condition for the existence of primitive solutions in
integers is the existence of solutions modulo $q$, and it is easy
to see that this is equivalent to ask for the existence of
solutions of $u^4 \equiv -1 \pmod{q}$, and by the classical theory
of cyclotomic fields we know that this congruence has solutions if
and only if $q \equiv 1 \pmod{8}$. Thus, for any odd prime $q$ not
verifying this condition (1.1) can not have primitive solutions,
for any value of $p$.
\\
From now on we will restrict to primes $q \equiv 1 \pmod{8}$ such
that $q \neq a^4 + b^4$ and we will give an algorithm that, using
the modularity of the Q-curve constructed by Darmon-Ellenberg
attached to a hypothetical solution of (1.1) and imitating the
Frey-Ribet approach, ends up for some values of $q$ in an
``impossible congruence", thus proving the non-existence of
primitive solutions.

\section{From non-existence of trivial solutions to non-existence of
solutions}

Fix a value of $q$ as above (an ``interesting" value) and assume
that for some $p > 13$ there exists a primitive solution $(A,B,C)$
of (1.1). We are assuming in particular that $q$ is odd (because
$2 = 1^4 + 1^4$), and it is easy to see that $C$ must be odd, and
prime to $3$. By assumption, $C \neq 1$. $A$ and $B$ have
different parity, thus we can assume that $A$ is even. Following
the construction of Darmon and Ellenberg we can attach, as we did
in [D] for the equation $x^4 + y^4 = z^p$, two Q-curves to the
triple $(A,B,C)$. For simplicity we will just consider one of them
(see, however, the remark at the end of this section concerning
the case $q=41$), the curve:
$$ E_{(A,B)}: \quad \quad y^2 = x^3 + 2(1+i) A x^2 + (-B^2 + i A^2) x $$
It follows from the results in [E], that $E_{(A,B)}$ is a degree
$2$ Q-curve defined over $\Q(i)$ semistable outside $2$. Odd
primes of bad reduction are precisely the primes dividing $q C$.
Thus, the curve has good reduction at $3$, and therefore it is
known to be modular (cf. [E-S]). Modularity can be interpreted in
``Galois language"  as follows: consider the compatible family of
two-dimensional Galois representations $\rho_\lambda$ of the full
Galois group of $\Q$ attached to (the Weil restriction of)
$E_{(A,B)}$. There exists a weight $2$ cuspidal modular form $f$
such that the representations $\rho_\lambda$ are attached to it,
i.e., the trace $a_t$ of the image of $\Frob \; t$ agrees with the
$t$-th Hecke eigenvalue of $f$, for every prime $t \nmid 2 q C$.
We know that the representations $\rho_\lambda$ have coefficients
in $\Q(\sqrt{2})$, and moreover that their restrictions to the
Galois group of $\Q(i)$ have rational coefficients. Thus, $f$ has
an inner twist, i.e., $f^{\sigma} = f \otimes \phi $, where
$f^{\sigma}$ is the Galois conjugate of $f$, and $\phi$ is the
quadratic character
corresponding to $\Q(i)$.\\
We list some properties of the representations $\rho_\lambda$ that
are proved in [E]: The $2$-part of the conductor of this family is
$32$ (recall that $A$ is even), the residual representations
$\bar{\rho}_\lambda$, for every $\ell > 13$, $\lambda \mid \ell$,
are irreducible, and their (projective) images can not fall in the
normalizer of a split Cartan subgroup of $\PGL_2(\mathbb{F}_\ell)$
(these two properties generalize similar results proved by Mazur
and Momose for
 elliptic curves over $\Q$).\\
 \newline
The close relation between the discriminant of $E_{(A,B)}$ and $q
C^p$ is a key point, allowing us to apply Frey's trick: Every odd
prime $q' \neq q$ of bad reduction of the curve disappears when
considering the (irreducible) residual representation
$\bar{\rho}_P$, $P \mid p$: this residual representation has
conductor $32 q$. Then, applying Ribet's level-lowering result, we
conclude that there exists a weight $2$  newform $f_2$ of level
$32 q$ such that, if $\{ b_t \}$ is the set of Hecke eigenvalues
of $f_2$ and $\{a_t \}$ the one of the modular form $f$ associated
to $E_{ (A,B) }$, we have the congruence:
 $$a_t \equiv b_t
\pmod{P} \qquad \qquad (2.1)$$
 for every prime $t \nmid 2 q C$.\\
 As $(A,B,C)$ is a
primitive solution of (1.1), every prime dividing $C$ is of the
form $8n+1$; in particular, the level of $f$ does not contain any
prime $t \equiv 3 \pmod{4}$.\\
\newline
Now, we want to derive a contradiction from the above $\mod \; P$
congruences by just looking at a few primes $t \equiv 3 \pmod{4}$.
On the one hand, we can control the values of the eigenvalues
$a_t$ of the modular form $f$ (whose level we do not know) for
such values of $t$ because $f$ has an inner twist given by the
$\mod \; 4$ character, hence these numbers must verify $a_t = c
\sqrt{2}$, for a rational integer $c$ that we can easily bound
using $|a_t| \leq 2 \sqrt{t}$. For example, we know that $a_3$
must be one of the following: $0, \pm \sqrt{2}, \pm 2 \sqrt{2}$.
On the other hand $f_2$ is new of level $32 q$, so by looking at
the tables in [S] (or by applying the Hecke software developed by
W. Stein in Magma) we can compute for
this finite set of eigenforms the first eigenvalues $b_t$.\\
The desired contradiction will be obtained as long as in the space
of modular forms of level $32 q$ no newform has its eigenvalues
$b_t$ for the first values of $t \equiv 3 \pmod{4}$ of the form $c
\sqrt{2}$. If this happens to be the case, then  congruences (2.1)
will only have solutions for a few primes $p$, and by
computing several $b_t$  one can shrink this set of primes.\\
\newline
Before computing some examples, we must stress that for certain
values of $q$ we know a priori that the method proposed will not
be successful. In fact, for a prime $q$ that can be written as a
sum of two biquadrates $q= A^4 + B^4$, the method is useless
because there will be a newform $f_2$ of level $32 q$ whose
eigenvalues $b_t$ behave precisely as above, namely the newform
attached to the Q-curve corresponding to the trivial solution
$(A,B,1)$. In fact, this newform will have field of coefficients
$\Q(\sqrt{2})$ and an inner twist by the character corresponding
to $\Q(i)$. Thus, a non-trivial solution $(A', B', C)$ of (1.1)
for such values of $q$, that would correspond to a newform $f$ of
higher level, may not exist, but we can not
answer this question with our method.\\
\newline
For a similar reason, our method does not work for primes $q$ that
can be written as $(2A)^4 + B^2$ (for example, $41 = 2^4 + 5^2$),
because even in this case a Q-curve can be attached (cf. [E]) to
the triple $(2A,B,1)$ and the corresponding newform again will
have level $32 q$ and the same  field of coefficients and inner
twist as above. However, there is still an opportunity to rescue
some of these primes: one can consider the second Q-curve attached
to a hypothetical solution $(A', B', C')$ of $\eq$, the curve
$E_{(B',A')}$ defined in [D] (we are assuming that $B'$ is odd).
The same method as above can now be applied, this time the
congruences relate the eigenvalues of the newform attached to this
Q-curve with those of a newform of level $256 q$, and one can hope
that in level $256 q$ there is no newform with coefficients in
$\Q(\sqrt{2})$ and an inner twist by the $\mod \; 4$ character
(the equation $A^4 + (2B)^2 = q$ has no solutions, thus there is
no obvious way to construct a Q-curve such that
 the corresponding newform has level $256 q$). This alternative
 method has not been tested by us in any example because it
 involves computations with too large levels: $256 q \geq 256 \cdot 41 > 10000
 $.\\

 \section{Three successful examples}

Let us consider the first three examples of ``interesting" primes
that can not be written as $(2A)^4 + B^2$: $$q=73, 89, 113$$
 We
have applied the algorithm described in the previous section to
solve the diophantine equations (1.1) for these values of $q$ and
$p>13$. To do this, we extracted from the tables in [S] the
eigenvalues $b_3 , b_7, b_{11}$ and $b_{19}$ for all newforms of
level $32 q$ for these values of $q$. For all the considered
newforms at least one of these eigenvalues is not of the form $c
\sqrt{2}$, $c$ an integer. Moreover, solving the congruences $a_t
\equiv b_t \pmod{P}$ where $a_t = c \sqrt{2}$, for every integer
$c$ such that  $ |a_t| = |c| \sqrt{2} \leq 2 \sqrt{t}$ and $t =
3,7,11$ and $19$, we obtain only one case where these congruences
have a solution: $q= 73$, $p = 17$ and the newform of level $73
\cdot 32$ named as $2336L$ in [S] (to solve the congruences, we
just compute the primes
in the resultant between the minimal polynomials of $a_t$ and $b_t$). \\
Thus, except for the case $q=73$ and $p=17$, we obtain a
contradiction, which proves the following:

\begin{teo}
\label{teo:ejem} The diophantine equations $\eq$ do not have
primitive solutions for $q= 73, 89$ and $113$ and $p>13$.
\end{teo}

Remark: The method proposed in this article can be used to compute
more examples. One can even ask the following question: Is there
an interesting prime $q$ (not a sum of two biquadrates) such that
the equation $x^4 + y^4 = q z^p$ has a primitive solution for some
$p>13$? What if we restrict to primes $q$ not of the form $A^4 +
B^2$?\\

 It remains to prove the theorem for $q=73$ and $p=17$. We will see that the
 non-existence of solutions in this case follows from the results
 in [E].\\
 What happens to the $\mod \; 17$ representation
 attached to the newform $2336L$ is that it has image contained in
 the normalizer of a Cartan subgroup (in fact, for this modular form congruences (2.1) hold modulo $17$ precisely when
  we
  take $a_t = 0$, for $t \equiv 3 \pmod{4}$). Moreover, this is due to a
 congruence between $2336L$ and the newform of level $32$. The
 congruence between these two cusp forms can be proved by direct
 computation, checking that it holds for sufficiently many
 eigenvalues  and applying Sturm's bound. The image is
 dihedral because of this congruence, since the newform of level
 $32$ has Complex Multiplication: the Galois representations attached to it are reducible
 when restricted to the Galois group of $\Q(i)$. As $17 \equiv 1
 \pmod{4}$, we see that the  $\mod \; 17$ representation
 attached to the newform of level $32$  has its image contained in
 the normalizer of a split Cartan subgroup. Now consider the
 newform $f$ corresponding to a solution of (1.1) with $q=73$ and $p=17$. If
 we suppose that we have a congruence $\mod \; 17$ between $f$ and the newform
 $2336L$, then we obtain a congruence also between $f$ and the
 newform of level $32$, and in particular we conclude that the
 image of the $\mod \; 17$ representation attached to $f$ will be
 contained in the normalizer of a split Cartan subgroup, but this
 contradicts Ellenberg's generalization to Q-curves of results of
 Momose (see section 2).\\
 \newline
  Another way of obtaining a contradiction
 without applying the results of Ellenberg-Momose is the
 following: the $\mod \; 17$ representation attached to the
 Q-curve can not have conductor $32$, its conductor must be $32 \cdot 73$, because it has semistable reduction at primes
  above $73$
 and the formula for the
 discriminant of this Q-curve shows that the prime $17$ does not
 divide the number of connected components of the reduction of the
 curve modulo a prime of $\Q(i)$ dividing $73$ (the exact power of such a prime dividing the discriminant is not
 a multiple of $17$).

\section{Bibliography}

[D] Dieulefait, L., {\it Modular congruences, Q-curves, and the diophantine equation $x^4 + y^4 = z^p $
}, preprint; available at:\\
http://www.math.leidenuniv.nl/gtem/view.php (preprint 55)\\

[E] Ellenberg, J.,{\it Galois representations attached to Q-curves
and the generalized Fermat equation $A^4 + B^2 = C^p$}, preprint; available at:\\
http://www.math.princeton.edu/~ellenber/papers.html\\

[ES] Ellenberg, J., Skinner, C., {\it On the modularity of
Q-curves},
Duke Math. J., 109 (2001), 97-122\\

[S] Stein, W., {\it The Modular Forms Database}, available at: \\
http://modular.fas.harvard.edu/Tables/index.html

%[S2] Stein, W., {\it },

\end{document}